\documentclass{article}[12pt]

\usepackage[english]{babel}
\usepackage[latin1]{inputenc}
\usepackage[T1]{fontenc}

\usepackage{amsmath, amsthm}
\usepackage{amscd}
\usepackage{amssymb}

\usepackage{url}

\usepackage{color}

\usepackage{enumerate}
\usepackage{graphicx}
\usepackage{float}

\usepackage{lipsum}

\newtheorem{thm}{Theorem}[section]
\newtheorem{theorem}[thm]{Theorem}
\newtheorem{claim}[thm]{Claim}

\newtheorem{corollary}[thm]{Corollary}
\newtheorem{lemma}[thm]{Lemma}
\newtheorem{proposition}[thm]{Proposition}
\newtheorem{definition}[thm]{Definition}

\theoremstyle{definition}
\newtheorem{example}[thm]{Example}
\newtheorem{examples}[thm]{Examples}

\newtheorem{remark}[thm]{Remark}

\begin{document}

\centerline{\Large \bf Superalgebraic methods in the classical theory}

\bigskip

\centerline{\Large \bf of representations. }

\bigskip

\centerline{\Large \bf  Capelli's identity, the Koszul map}

\bigskip

\centerline{\Large \bf  and the center of the enveloping algebra $\mathbf{U}(gl(n)).$}

\bigskip

\centerline{A. Brini}
\centerline{\it $^\flat$ Dipartimento di Matematica, Universit\`{a} di
Bologna }
 \centerline{\it Piazza di Porta S. Donato, 5. 40126 Bologna. Italy.}
\centerline{\footnotesize e-mail: andrea.brini@unibo.it}
\medskip

\section{Introduction}

The theme of the superalgebraic extension of many theories both in Geometry
and Algebra has noble origins that lay in Physics and has been at the core of a wide range of works
by several prominent mathematicians during the last three decades (see, e.g.  \cite{FSS-BR}).
We are not concerned with the general theory here, but
we limit ourselves to show how the use of superalgebraic methods sheds new light on some
classical themes of representation theory and it leads to
significant simplifications of  traditional proofs.

\medskip

In this paper we essentially deal with Capelli's identity\footnote{ The celebrated identity that now bears his name was proved
in $1887$  \cite{Cap1-BR}, when he held the chair of Algebra at
the University of Naples. Alfredo Capelli graduated from the University of Rome in 1877 and
then continued to develop his mathematical skills working as
Felice Casorati's assistant at the University of Pavia.
Casorati corresponded regularly with Weierstrass and had developed a strong link between Italian and German mathematicians.
Capelli spent time at the University of Berlin where he was influenced by Karl Weierstrass and Leopold Kronecker.}
and with the center of the enveloping algebra of the
general linear Lie algebra $gl(n)$.

\medskip

Capelli's identity is the cornerstone of
the classical theory of algebraic invariants  (see, e.g \cite{Weyl-BR}, \cite{Procesi-BR}).
Distinguished contemporary authors referred to this identity as ``mysterious'' (see, e.g.  \cite{ABP-BR},  \cite{Howe-BR}) , and it still
provides a quite active research field.
The main mathematical object in  Capelli's identity is a ``deteminantal'' operator that results (in modern language) in a central
element of the enveloping algebra $\mathbf{U}(gl(n))$ of the general Lie algebra $gl(n).$ The generalization of this operator
led Capelli to study the mathematical
structure that we nowadays call the {\it{center} of } $\mathbf{U}(gl(n))$, to prove that it is a polynomial algebra and
to explicitly describing a family of free algebraic generators (\cite{Cap2-BR}, \cite{Cap3-BR}, $1893$).

\medskip

By way of elementary motivation, we  show that the proof of  Capelli's identity can be
reduced to a straightforward computation just by using a touch of superalgebraic notions.

This point of view is extended to the study of the enveloping algebra $\mathbf{U}(gl(n))$. The notions of {\it{determinantal and
permanental Capelli bitableaux}}
provide two relevant classes of bases that arise from  {\it{Straightening Laws}} (see. e.g. \cite{Brini3-BR}, \cite{drs-BR}, \cite{rota-BR},
\cite{Procesi-BR}).

In the final section, we submit new results on the center of $\mathbf{U}(gl(n))$.
These results - which cannot even be expressed
without appealing to the superalgebraic notation  - allows a variety of classical and recent
results to be almost trivially proved and  be put under one  roof.

\medskip

The crucial methodological tool of our approach to Capelli's identity and to the center of the enveloping algebra $\mathbf{U}(gl(n))$
is the {\it{superalgebraic method of virtual variables}}, which is, in turn, an extension
of Capelli's method of  {\it{ variabili ausilarie}}. Capelli introduced the method of \
``variabili ausiliarie'' in order
to manage symmetrizer operators in
terms of polarization operators, to  simplify the study of some skew-symmetrizer operators (namely, the famous central Capelli operator) and
developed this idea in a systematic way in his beautiful treatise \cite{Cap4-BR}.
As a matter of fact, Hermann Weyl's  sketchy proof of  Capelli's identity  (\cite{Weyl-BR}, page $39$ ff.) must be regarded as
an application of the method of \  ``variabili ausiliarie'' in disguise (we refer the reader to the footnote to Theorem \ref{Capelli identity} for further details).

Unfortunately, Capelli's idea was well suited to
treat symmetrization, but it did not work in the same efficient way while dealing with skew-symmetrization.

We had to wait the
 introduction of the notion of {\it{superalgebras}} to have the right conceptual framework to treat symmetry and skew-symmetry in one and the same way.
To the best of our knowledge, the first mathematician who intuited the connection between Capelli's idea and superalgebras was
Koszul in $1981$  \cite{Koszul-BR}; Koszul
 proved that the classical determinantal Capelli operator can be rewritten - in a much simpler way - by adding to the symbols to be dealt with
an extra auxiliary symbol that obeys to different commutation relations.

The {\it{supersymmetric method of virtual variables}} was developed in its full extent and generality (for the   general linear
Lie superalgebras $gl(m|n)$ - in the notation
of \cite{KAC1-BR}) in the series of notes \cite{Brini1-BR}, \cite{Brini2-BR}, \cite{Brini3-BR}, \cite{Brini4-BR} by Teolis and the present author.

\medskip

This note is organized as follows. In section $2$, we summarize some
elementary material required
to define the classical Capelli operator  and to state Capelli's identity.

In section $3$, we introduce the minimum of superalgebraic concepts and notations required
to express the Capelli operator  in a compact way and to provide a few lines proof of classical Capelli's identity.

In Section $4$, we provide a systematic treatment of the
method of virtual variables.
The starting point of the method is to  introduce new symbols, called {\it{ virtual variables}}
(``variabili ausiliarie'' in the language of Capelli \cite{Cap4-BR}),
which may have different signature
(parity, $\mathbb{Z}_2-$degree) than the signature  of the
symbols to be dealt with. This leads to a Lie superalgebra $gl(m_1|m_2+n)$, where $m_1$ is the number of positive virtual symbols and
 $m_2$ is the number of negative virtual symbols (informally, we assume $m_1, m_2$ "sufficiently large").

The main technical device of the method is the notion of {\it{virtual
algebra}} $Virt(m_1+m_2,n)$
as a subalgebra of the enveloping algebra ${\mathbf{U}}(gl(m_1|m_2+n))$.
The virtual algebra $Virt(m_1+m_2,n)$ provides an effective method to define special elements  and to derive
deep identities in ${\mathbf{U}}(gl(n))$.

Specifically, ${\mathbf{U}}(gl(n))$ is the image of
$Virt(m_1+m_2,n)$ under an algebra homomorphism - the Capelli epimorphism -
whose kernel is the linear span of a set of
monomials characterized by a simple combinatorial property;
Capelli rows and bitableaux, while having no simple expressions
in ${\mathbf{U}}(gl(n))$, admit very simple preimages in $Virt(m_1+m_2,n)$
and, hence, all the calculations are carried on in the virtual
algebra, modulo the kernel of the Capelli epimorphism.

In section $5$, the sets of {\it{determinantal}} and {\it{permanental}} {\it{Capelli bitableaux}} are introduced;
these elements of $\textbf{U}(gl(n))$ could not be defined, nor imagined without recourse to the method of virtual variables.

Capelli bitableaux give rise to {\it{straightening laws}} for the enveloping algebra \cite{Brini3-BR} that are in all respect similar to
the ordinary ones for the classical (symmetric) algebra of algebraic forms ${\mathbb C}[M_{n,n}] \cong Sym[gl(n)]$
(see, e.g. \cite{drs-BR}, \cite{DKR-BR}, \cite{DEP-BR}, \cite{rota-BR}).
We developed this connection between  $Sym[gl(n)]$ and  $\textbf{U}(gl(n))$ by introducing a linear isomorphism
$\mathfrak{T}: {\mathbb C}[M_{n,n}] \cong Sym[gl(n)] \rightarrow \textbf{U}(gl(n))$,  called the {\it{bitableau correspondence}}  \cite{Brini3-BR};
the isomorphism $\mathfrak{T}$  maps each bitableau
of ${\mathbb C}[M_{n,n}]$ to the Capelli bitableaux of $\textbf{U}(gl(n))$
parametrized by the same pair of Young tableaux, both in the determinantal and in the permanental cases. The sets of semistandard determinantal
and co-semistandard permanental bitableaux yield two remarkable bases of $\textbf{U}(gl(n))$; furthermore, the map $\mathfrak{T}$ turns out to be the inverse of the
{\it{Koszul map}} $\mathfrak{K}$ (\cite{Koszul-BR}, \cite{Brini4-BR}).

The final section is devoted to the study of the center $\mathfrak{Z}(\textbf{U}(gl(n)))$ of the enveloping algebra $\textbf{U}(gl(n)).$
The  new concept of {\it{rectangular Capelli-Deruyts tableau}}  is introduced and a new result  (Theorem \ref{Teorema di sfilamento-BR})
is presented.

Due to their virtual presentation, rectangular Capelli-Deruyts tableaux are almost immediately recognized to be central elements
in  $\textbf{U}(gl(n))$ and,
by Theorem \ref{Teorema di sfilamento-BR}, they can be expanded in terms of the central generators first discovered by Capelli
in $1893$ (\cite{Cap2-BR}, \cite{Cap3-BR}).

By combining this result with the beautiful description of the Harish-Chandra isomorphism in the terms of {\it{shifted symmetric polynomials}}
(Okounkov\footnote{Fields Medalist 2006} and Olshanski  \cite{OkOlsh-BR}),  a variety of  results
(see. e.g. \cite{Brundan-BR}, \cite{Cap2-BR}, \cite{Cap3-BR}, \cite{Howe-BR},
\cite{HU-BR}, \cite{Molev-BR}, \cite{OkOlsh-BR}, \cite{Procesi-BR}, \cite{UmedaCent-BR},  \cite{Umeda-BR}) on free sets of generators of the polynomial algebra
$\mathfrak{Z}(\textbf{U}(gl(n)))$ is deduced.

\medskip

We extend our heartfelt thanks to Rita Fioresi, Alberto Parmeggiani, Francesco Regonati and Antonio G.B. Teolis for their encouragement, advice and invaluable suggestions.

\section{
%Polarizations and Capelli's operator}
The classical Capelli identity}

%\comment{R: ho fatto def}
In this section, we   essentially refer to \cite{Procesi-BR}.

The \textit{algebra of algebraic forms in $n$ vector variables
of dimension $d$} is the polynomial algebra in $n \times d$ variables:
$$
{\mathbb C}[M_{n,d}] :=    {\mathbb C}[x_{ij}]_{i=1,\ldots,n; j=1,\ldots,d}
$$
where $M_{n,d}$ represents the  matrix
with $n$ rows and $d$ columns with ``generic" entries $x_{ij}$:
$$
M_{n,d} = \left[ x_{ij} \right]_{i=1,\ldots,n; j=1,\ldots,d}=
 \left[
 \begin{array}{ccc}
 x_{11} & \ldots & x_{1d} \\
 x_{21} & \ldots & x_{2d} \\
 \vdots  &        & \vdots \\
 x_{n1} & \ldots & x_{nd} \\
 \end{array}
 \right].
$$
%\end{definition}

We can interpret $x_{ij}$, with a slight abuse of
notation, as the {\it{multiplication operator}}:
$$
x_{ij} : {\mathbb C}[M_{n,d}] \rightarrow {\mathbb C}[M_{n,d}],
\quad x_{ij}(f) = x_{ij} \cdot f, \ \forall f \in {\mathbb C}[M_{n,d}].
$$

The \textit{Weyl algebra}
$$
\mathbf{W}_{n,d}:=  {\mathbb C}[x_{ij}, \frac {\partial} {\partial
x_{ij}}] \subseteq \mathrm{End}_{{\mathbb C}}[{\mathbb C}[M_{n,d}]]
%\langle x_i|j, \partial_{x_h|k} \rangle
$$
is the subalgebra of $\mathrm{End}_{{\mathbb C}}[{\mathbb C}[M_{n,d}]]$ generated by
the multiplication operators $x_{ij}$ and the partial derivative operators $\frac {\partial} {\partial x_{ij}}.$

This is not a commutative algebra:
$$
 \left[ x_{ij}, \frac {\partial} {\partial x_{hk}} \right] = \delta_{ih}\delta_{jk} \texttt{I},
$$
where $\texttt{I}$ denotes the identity operator.

In $\mathbf{W}_{n,d}$, we consider the \textit{polarization operator}:
$$
D_{x_h,x_k}=\sum_{j=1}^d  x_{hj} \frac {\partial} {\partial
x_{kj}}.
%\sum x_h|j \partial_{x_k|j}
$$

The operator $D_{x_h,x_k}$ is characterized as the unique derivation of the algebra
${\mathbb C}[M_{n,d}]$ such that
$$
D_{x_{h},x_{k}} \left( x_{ij} \right) = \delta_{ki} \ x_{hj}  \qquad
\forall j = 1, \ldots, d.
$$

The  subalgebra (with $1$) of $ \textbf{W}_{n,d}$  generated by the polarization operators is
denoted by $ \textbf{P}_{n,d}$ and called the
\textit{algebra of polarizations}.

The following identity holds:
$$
[D_{x_i,x_j},D_{x_h,x_k}] = \delta_{jh}D_{x_i,x_k}-\delta_{ik}D_{x_h,x_j}
$$

One may recognize the commutation relations of elementary matrices, so
we have a representation $\rho:\mathbf{U}(gl(n)) \longrightarrow \mathbf{W}_{n,d}$. If $n \leq d$
we have an injection, so we obtain an isomorphism:
$$
\rho(\mathbf{U}(gl(n))\cong \mathbf{P}_{n,d}
$$
%where $P_{n,d}$ is the subalgebra of  the polarization operators.

Given a square matrix
$
\mathbf{M_n} = \left[ a_{ij} \right]_{i, j =1,\ldots,n}
$
with entries $a_{ij}$ in a noncommutative algebra, we will consider its {\it {column determinant}}
$$
\textbf{cdet}(\mathbf{M_n}) = \sum_{\sigma \in \mathcal{S}_n} \ (-1)^{| \sigma |} \
a_{{\sigma(1)},{1}}a_{{\sigma(2)},{2}} \cdots a_{{\sigma(n)},{n}}.
$$

The \textit{classical Capelli operator} $H_{n,d}$ is the element in
$\mathbf{P}_{n,d} \subset \mathbf{W}_{n,d}$ given by the column determinant:

\fontsize{10}{12}\selectfont

$$
H_{n,d} = \textbf{cdet}
 \left(
 \begin{array}{cccc}
 D_{x_{1},x_{1}} + (n-1)\texttt{I} & D_{x_{1},x_{2}} & \ldots &  D_{x_{1},x_{n}} \\
 D_{x_{2},x_{1}} & D_{x_{2},x_{2}} + (n-2)\texttt{I} & \ldots  &  D_{x_{2},x_{n}}\\
 \vdots  &                                &  &\vdots \\
 %D_{x_{n-1},x_{1}} & D_{x_{n-1},x_{2}} & \ldots &   D_{x_{n-1},x_{n}}\\
D_{x_{n},x_{1}} & D_{x_{n},x_{2}} & \ldots & D_{x_{n},x_{n}}\\
 \end{array}
 \right) =
$$
$$
= \textbf{cdet} \left( \ D_{x_{i},x_{j}} + \delta_{ij}(n-i)\texttt{I} \ \right)_{i, j = 1, 2, \ldots, n}.
$$

The operator $H_{n,d}$ determines a unique element $\mathbf{H}_n$ in $\mathbf{U}(gl(n))$
(with $n \leq d$) which is central (see, e.g.  \cite{Umeda-BR}).

The
\textit{Cayley operator} $\Omega_n$ (see, e.g. \cite{TURN-BR},  \cite{Weyl-BR}):
$$
\Omega_n =
 \textbf{det}
 \left( \frac {\partial} {\partial x_{ij}} \right)_{i, j = 1, 2, \ldots, n}.
$$
Notice that $\Omega_n$ is computed as a ``true'' determinant,
since its entries are commuting operators.

The \textit{bracket}\footnote{This notation and terminology is due to Cayley.} is the element
\begin{equation}\label{bracket-BR}
[x_1, \ldots, x_n] = \mathrm{\bf{det}} \left[ x_{ij} \right]_{i, j =1,\ldots,n} \in {\mathbb C}[M_{n,d}]
\end{equation}
We are ready to state  the
{\it { Capelli's identity}}.
\begin{theorem}\label{Capelli identity}(Capelli, \cite{Cap1-BR} 1887)\label{capelli-thm-BR}
$$
H_{n, d}(f) = \begin{cases} 0 &\mbox{if } n > d \\
[x_1, \ldots, x_n] \Omega_n(f) & \mbox{if } n = d, \end{cases}
\quad \forall f \in {\mathbb C}[M_{n,d}].
$$
\end{theorem}
The case $n=d$ is called the \textit{special identity} (\cite{Weyl-BR}).
The special identity\footnote{As we announced in the Introduction, we formalize Weyl's  cumbersome
argument (\cite{Weyl-BR}, page 39)
in the special case $n = d = 2.$
Add to the vector symbols $x_1 = (x_{11}, x_{12})$ and $x_2 = (x_{21}, x_{22})$ two copies (variabili ausiliarie in Capelli's language)
$x^{'}_1 = (x^{'}_{11}, x^{'}_{12})$ and $x^{'}_2 = (x^{'}_{21}, x^{'}_{22})$, and embed the polynomial algebra
${\mathbb C}[M_{2,2}] = {\mathbb C}[(x_i|j)]_{i, j = 1,2}$
we are  concerned with into the enlarged polynomial algebra
${\mathbb C}[(x_i|j), (x^{'}_i|j)]_{i, j = 1,2}$ (these are the operators that Weyl denotes by the symbol $\Delta_{x_{i},x_{j}},$ and whose composition
Weyl called ``pseudo-composition''). The polarization operators $D_{x^{'}_{i},x_{j}}$ COMMUTES, and, therefore we may consider the ordinary
determinant $det[D_{x^{'}_{i},x_{j}}].$
When we apply the operator
$$
H^{'}_2 = D_{x_{1},x^{'}_{1}}D_{x_{2},x^{'}_{2}} \ det[D_{x^{'}_{i},x_{j}}]
$$
to an algebraic form $f \in {\mathbb C}[M_{2,2}]$,
 we get the form $[x_1, x_2] \ \Omega_2(f).$  By applying the commutator identities in the enlarged Weyl algebra $\mathbf{W}_{4,2}$,
we find that the column determinant
$
H_{2, 2} = \textbf{cdet} \left( \ D_{x_{i},x_{j}} + \delta_{ij}(2-i)\texttt{I} \ \right)_{i, j = 1, 2}
$
satisfies the identity $H_{2, 2} = H^{'}_2 + T$, where $T$ is an operator that acts trivially on the proper algebra
${\mathbb C}[M_{2,2}] = {\mathbb C}[(x_i|j)]_{i, j = 1,2}$. In plain words, we can subsitute the action of the column determinant
$H_{2, 2}$ - which lives in a noncommutative algebra - with the action of the true determinant (in commuting operators) $H^{'}_2.$} may be expressed as the following identity
in the Weyl algebra:
$$
H_{n, n} = [x_1, \ldots, x_n] \Omega_n.
$$

\section{The superalgebraic method of virtual variables for
 Capelli's identities}

For the sake of readability, from now on, we will write $(x_i|j)$ in place of $x_{ij}.$

Besides the (commuting) variables $(x_i|j), \ i = 1, \ldots, n, \ j =1, \ldots, d$, we consider a new series
of variables $(\alpha|j), \ j =1, \ldots, d$ and we assign
to them parity equal to one, while we assume the parity of $(x_i|j)$
to be zero:
$$
|(x_i|j)| = 0 \in \mathbb{Z}_2, \qquad |(\alpha|j)|  = 1 \in \mathbb{Z}_2.
$$
So our variables are now $\mathbb{Z}_2$-graded. The elements with parity
zero are also called \textit{even}, while those with parity one
are called \textit{odd}.

The symbol $\alpha$ is called \textit{virtual (or auxiliary) symbol}.

In analogy with the ordinary
setting that we have discussed in the previous section,
we define the \textit{supersymmetric algebra}:
$$
{\mathbb C}[M_{1|n,d}] := {\mathbb C}[(\alpha|j),
(x_i|j)]_{i=1,\ldots,n; j=1,\ldots,d}
$$
as the algebra generated by the $\mathbb{Z}_2$-graded variables $(\alpha|j),
(x_i|j)$
subject to the  commutation relations:
$$
(x_h|i)(x_k|j)  =  (x_k|j)(x_h|i), \ (x_h|i)(\alpha|j)  =  (\alpha|j)(x_h|i), \ (\alpha|i)(\alpha|j)  =  -(\alpha|j)(\alpha|i).
$$

Strictly speaking we have:
$$
{\mathbb C}[M_{1|n,d}] \cong  \bigwedge \left[ (\alpha|j) \right]
\bigotimes      {\mathrm{Sym}} \left[ x_{hj} \right],
%k[M_{m|n,d}]=\wedge(\alpha|j)k[(x_k|j)]
$$
and, therefore, ${\mathbb C}[M_{1|n,d}]$ is a $\mathbb{Z}_2-$graded algebra, whose
$\mathbb{Z}_2-$graduation is  inherited by the natural one in the exterior algebra.

Given $a, b \in \{x_1, \ldots, x_n, \alpha \}$, we can define, similarly to what we did before, the operator
$D_{a,b}$ which is a \textit{superpolarization}:
$$
D_{a,b}\left( (c|j) \right) = \delta_{b,c} \ (a|j),
\qquad \forall j = 1, \ldots, d.
$$
In general, $D_{a,b}$ is not a derivation, but a
\textit{superderivation} (see, e.g \cite{FSS-BR}, \cite{Scheu-BR}); in other words, we
need to introduce a sign in the Leibniz rule
if we want to obtain a consistent definition of the  operators $D_{a,b}$:
$$
D_{a,b}(FG)=D_{a,b}(F)G+ (-1)^{|D_{a,b}||F|}FD_{a,b}(G)
$$
where:
$
|F| \in \mathbb{Z}_2
$
denotes the $\mathbb{Z}_2-$degree (or, {\it{parity}}) of an $\mathbb{Z}_2-$homogeneous element $F \in {\mathbb C}[M_{1|n,d}]$ and
$$
|D_{{x_h},{x_k}}| = |D_{{\alpha},{\alpha}}| = 0 \in {\mathbb Z}_2, \quad
\qquad  |D_{{\alpha},{x_k}}| = |D_{{x_h},{\alpha}}| = 1 \in {\mathbb Z}_2.
$$

The space spanned by the superpolarizations form a Lie (sub)superalgebra of the Lie superalgebra
$\mathrm{End}(\mathbb C[M_{1|n}])$ (see, e.g. \cite{FSS-BR}, \cite{Scheu-BR}, \cite{KAC1-BR}), where the Lie (super)bracket is defined
introducing also a sign:
$$
[A,B]=AB-(-1)^{|A||B|}BA
$$
This is called the \textit{supercommutator}.

\begin{claim}

Given $a, b , c, d\in \{x_1, \ldots, x_n, \alpha \}$, we have
$$
[D_{a,b}, D_{c,d}] = \delta_{b,c} D_{a,d} + (-1)^{|D_{a,b}| |D_{c,d}|} \delta_{a,d} D_{c,b}.
$$

\end{claim}

There is a fundamental yet simple fact coming from a
straightforward calculation ($[x_1, \ldots,  x_n]$ denotes the {\it{bracket}},  eq. (\ref{bracket-BR}) ):
\begin{lemma}
\begin{equation}\label{virtual presentation of brackets-BR}
D_{x_n, \alpha} D_{x_{n - 1}, \alpha} \cdots D_{x_1, \alpha} \left(  (\alpha|1) \ldots (\alpha|n - 1)(\alpha|n) \ \right)  =
[x_1, \ldots,  x_n].
\end{equation}
\end{lemma}

Consider now the {\it{product}} of superpolarizations:
$$
\mathcal{H}_n = D_{x_n,\alpha} \dots D_{x_1,\alpha}D_{\alpha,x_1} \dots D_{\alpha,x_n}
\in \mathrm{End}(\mathbb C[M_{1|n}]).
$$
\begin{claim}
The subalgebra ${\mathbb C}[M_{n,d}] \hookrightarrow  {\mathbb C}[M_{1|n,d}]$
is invariant for the operator $\mathcal{H}_n$.
\end{claim}

Denote by $\mathcal{H}_n \mid_{ {\mathbb C}[M_{n,d}]}$ the restriction of the operator
$\mathcal{H}_n$ to the subalgebra ${\mathbb C}[M_{n,d}].$

\begin{claim}(see, e.g. \cite{Brini2-BR}, \cite{Bri-BR})

The restriction $\mathcal{H}_n \mid_{ {\mathbb C}[M_{n,d}]}$ can be
 expressed as an operator in the  polarizations
$D_{{x_h}, {x_k}},$  $\ h, k = 1, \ldots, n.$ In plain words,
it belongs to the classical algebra of polarizations
$\textbf{P}_{n,d}   \hookrightarrow \textbf{W}_{n,d}$.

\end{claim}

$\mathcal{H}_n$ is written as a monomial operator, but its restriction $\mathcal{H}_n \mid_{ {\mathbb C}[M_{n,d}]}$
becomes quite different when
expressed in terms of classical polarizations. This  fact has many consequences
in terms of computations, but let us first understand the close
connection between $\mathcal{H}_n$ and the classical Capelli operator.

\medskip

The next result is a special case of the result we called the  ``Laplace expansion for Capelli rows'' ( \cite{BRT-BR} Theorem $2$,
\cite{Bri-BR} Theorem $6.3$).

\begin{theorem}\footnote{A sketchy proof of this result can also be found in \cite{Koszul-BR}. From a representation theoretic point of view, Theorem \ref{Presentazione virt-BR} is
a special case of results that are discussed in the final section of the present paper, to which we refer the reader.}\label{Presentazione virt-BR}

For all $d \in \mathbb{Z}^+,$ we have:
$$
\mathcal{H}_n \mid_{ {\mathbb C}[M_{n,d}]} = \textit{H}_{n, d}.
$$
\end{theorem}

\medskip

The proof of Capelli's identities reduces to a
straightforward computation.

Let us first consider the action of the {\sl virtual} part
of $\mathcal{H}_n|_{\mathbb C[M_{n,d}]}$, namely:
$$
 D_{\alpha, x_1}  \cdots D_{\alpha, x_{n - 1}}D_{\alpha, x_n}
$$
on the generic monomial:
$$
\textbf{m}  = \prod_{i = 1}^n \ (\ \prod_{j = 1}^d  \ (x_i| j)^{d_{i j}} \ ) \in {\mathbb C}[M_{n,d}], \qquad d_{i j} \in \mathbb{N}.
$$
Notice that every monomial in the expression
\begin{equation} \label{monomial-BR}
 D_{\alpha, x_1}  \cdots D_{\alpha, x_{n - 1}}D_{\alpha, x_n} (\textbf{m})
\end{equation}
contains exactly $n$ occurrences of the anticommuting variables of
type:
$$
(\alpha|j), \qquad j = 1, \ldots, d.
$$

\medskip\noindent
\textbf{Case $n > d.$}
To prove the first identity it is enough to
observe that all the monomials resulting in
the expression (\ref{monomial-BR}) contain
the square of an odd variable, hence they are all zero. So
$$
\mathcal{H}_n(\textbf{m}) = \textit{H}_{n, d} (\textbf{m}) = 0.
$$
Hence we have the first identity, by linearity.

\medskip

\medskip\noindent
\textbf{Case $n = d.$}
Let us set
$(x_i|j)^{d_{i j} -1} =
0$ if $d_{i j} = 0$, and recall $(\alpha|j)^2 = 0$
for $j = 1, 2, \ldots,d.$
We have:
\begin{multline*}
D_{\alpha, x_1}  \cdots D_{\alpha, x_{n - 1}}D_{\alpha, x_n} (\textbf{m}) =  \\ =
D_{\alpha, x_1}  \cdots D_{\alpha, x_{n - 1}}D_{\alpha, x_n} \left( \prod_{i = 1}^n \
(\ \prod_{j = 1}^d  \ (x_i| j)^{d_{i j}} \ ) \right)      \\ =
\sum_{\sigma \in {\mathcal{S}}_n} \ \left( \prod_{i, j = 1}^n \ d_{i,\sigma_i}
(x_i|j)^{d_{i  j}- \delta_{j, \sigma_i}} (\alpha|\sigma_1) \cdots
(\alpha|\sigma_{n-1}) (\alpha|\sigma_n) \right) =
\end{multline*}

\begin{multline*}
 = \left( \sum_{\sigma \in {\mathcal{S}}_n} (-1)^{|\sigma|}  \prod_{i, j = 1}^n \ d_{i,\sigma_i}
(x_i|j)^{d_{i  j}- \delta_{j, \sigma_i}} \right) \cdot (\alpha|1) \cdots (\alpha|n - 1)(\alpha|n)\\
\\ =
\Omega_n (\textbf{m}) \cdot (\alpha|1) \ldots (\alpha|n - 1)(\alpha|n) \\
\\ =
(\alpha|1) \ldots (\alpha|n - 1)(\alpha|n) \cdot \Omega_n (\textbf{m}).
\end{multline*}

We proved:
$$
D_{\alpha, x_1}  \cdots D_{\alpha, x_n} (\textbf{m}) =
(\alpha|1) \ldots (\alpha|n) \cdot \Omega_n (\textbf{m}).
$$

Applying the operator
$$
D_{x_n, \alpha} D_{x_{n - 1}, \alpha} \cdots D_{x_1, \alpha}
$$
to both sides we have:
\begin{multline*}
\textit{H}_{n,n } (\textbf{m}) = \mathcal{H}_n (\textbf{m}) =
D_{x_n, \alpha}  \cdots D_{x_1, \alpha} \cdot D_{\alpha, x_1}  \cdots D_{\alpha, x_n} (\textbf{m})  \\
\\ =
D_{x_n, \alpha} D_{x_{n - 1}, \alpha} \cdots D_{x_1, \alpha} \left( \ (\alpha|1) \ldots (\alpha|n - 1)(\alpha|n)
\cdot \Omega_n (\textbf{m}) \ \right)  \\
\\ =
\left(D_{x_n, \alpha} D_{x_{n - 1}, \alpha} \cdots D_{x_1, \alpha} \left( \ (\alpha|1) \ldots (\alpha|n - 1)(\alpha|n) \ \right) \right)
\cdot \Omega_n (\textbf{m})   \\
\\ =
[x_1, \ldots, x_{n - 1}, x_n] \cdot \Omega_n (\textbf{m}),
\end{multline*}
Hence by linearity we obtain the special identity.
\fbox

\section{The method of
 virtual supersymmetric variables for $\textbf{U}(gl(n))$}

 Let us consider the
vector spaces $V_n$ and the {\it{ auxiliary}} vector spaces
$V_{m_1}$ and $V_{m_2}$ (informally, we assume that
$dim(V_{m_1})=m_1$ and $dim(V_{m_2})=m_2$ are ``sufficiently large").
$V_n$ is called the space of {\it{proper vectors}}, and the spaces
$V_{m_1}$ and $V_{m_2}$ are called the spaces of {\it{even virtual
vectors}} and of  {\it{odd virtual vectors}}, respectively.

Let $W = W_0 \oplus
W_1$ be the $\mathbb{Z}_2-$graded vector space, where
$$
W_0 = V_{m_1}, \qquad W_1 = V_{m_2} \oplus V_n
$$
and let $gl(m_1|m_2+n)$ denote the general linear Lie superalgebra of $W = W_0 \oplus W_1$
(see, e.g. \cite{KAC1-BR}, \cite{Scheu-BR}, \cite{FSS-BR}).

Let
$
A_0 = \{ \alpha_1, \ldots, \alpha_{m_1} \},$  $A_1 = \{ \beta_1, \ldots, \beta_{m_2} \},$
$L = \{ x_1, \ldots, x_n \}$
denote distinguished  bases of $V_{m_1}$, $V_{m_2}$ and $V_n$, respectively; therefore $|\alpha_s| = 0 \in \mathbb{Z}_2,$
and $|\beta_t| = |x_i|   = 1 \in \mathbb{Z}_2.$

Let
$$
\{ e_{a, b}; a, b \in A_0 \cup A_1 \cup L \}, \qquad |e_{a, b}| =
|a|+|b| \in \mathbb{Z}_2
$$
be the standard $\mathbb{Z}_2-$homogeneous basis of $gl(m_1|m_2+n)$ provided by the
elementary matrices.

The supercommutator of $gl(m_1|m_2+n)$ has the following explicit form:
$$
\left[ e_{a, b}, e_{c, d} \right] = \delta_{bc} \ e_{a, d} - (-1)^{(|a|+|b|)(|c|+|d|)} \delta_{ad}  \ e_{c, b},
$$
$a, b, c, d \in A_0 \cup A_1 \cup L.$

\medskip

In analogy with the ordinary
setting  we  discussed in the previous section,
the \textit{supersymmetric algebra}
$$
{\mathbb C}[M_{m_1|m_2+n,d}]
$$
is the algebra generated by the ($\mathbb{Z}_2$-graded) variables $(\alpha_s|j), (\beta_t|j), (x_i|j)$,
 where
 $$
 |(\alpha_s|j)| = 1 \in \mathbb{Z}_2 \ \ and \ \  |(\beta_t|j)| = |(x_i|j)| = 0 \in \mathbb{Z}_2,
 $$

subject to the commutation relations:

\
$$
(a|h)(b|k) = (-1)^{|(a|h)||(b|k)|} \ (b|k)(a|h),
$$
for $a, b \in  \{ \alpha_1, \ldots, \alpha_{m_1} \} \cup \{ \beta_1, \ldots, \beta_{m_2} \} \cup \{x_1, x_2, \ldots , x_n\}.$

We have:
$$
{\mathbb C}[M_{m_1|m_2+n,d}] \cong  \bigwedge \left[ (\alpha_s|j) \right]
\bigotimes      {\mathrm{Sym}} \left[ (\beta_t|j), (x_h|j) \right],
%k[M_{m|n,d}]=\wedge(\alpha_s|j)k[(x_k|j)]
$$
and, therefore, ${\mathbb C}[M_{m_1|m_2+n,d}]$ is a $\mathbb{Z}_2-$graded algebra (superalgebra), whose
$\mathbb{Z}_2-$graduation is  inherited by the natural one in the exterior algebra.

Given two symbols $a, b \in A_0 \cup A_1 \cup L$, the {\it{superpolarization}} $D_{a,b}$ of $b$ to $a$
is the unique superderivation (see, e.g. \cite{FSS-BR}, \cite{Scheu-BR}, \cite{KAC1-BR})
 of ${\mathbb C}[M_{m_1|m_2+n,d}]$ of parity $|D_{a,b}| = |a| + |b| \in \mathbb{Z}_2$ such that
\begin{equation}
D_{a,b} \left( (c|j) \right) = \delta_{bc} \ (a|j), \ c \in A_0 \cup A_1 \cup L, \ j = 1, \ldots, d.
\end{equation}

We have a representation of $gl(m_1|m_2+n)$ on
${\mathbb{C}}[\mathbb{C}[M_{m_1|m_2+n,d}]]$ sending the elementary
matrices in the corresponding superpolarizations:
$$
\varrho : gl(m_1|m_2+n) \rightarrow End_\mathbb{C}[\mathbb{C}[M_{m_1|m_2+n,d}]]
$$
$$
e_{a,b} \rightarrow D_{a,b}, \qquad a, b \in A_0 \cup A_1 \cup L.
$$
Hence this defines a morphism (i.e. a representation):
$$
\varrho : \textbf{U}(gl(m_1|m_2+n)) \rightarrow End_\mathbb{C}[\mathbb{C}[M_{m_1|m_2+n,d}]].
$$

\begin{definition}{\bf{(Irregular expressions and the ideal $\mathbf{Irr}$)}}\label{Irregular expressions-BR}

We say that a product
$$
e_{a_mb_m} \cdots e_{a_1b_1} \in \textbf{U}(gl(m_1|m_2+n))
$$
is an {\it{irregular expression}} whenever
 the following condition on the occurrences of the virtual
variables $\alpha_s$ and $\beta_t$ holds: there exists a right subsequence
$e_{a_i,b_i} \cdots e_{a_2,b_2} e_{a_1,b_1}$, $i \leq m$ and a
virtual symbol $\gamma \in A_0 \cup A_1$ such that
\begin{equation}
 \# \{j;  b_j = \gamma, j \leq i \}  >  \# \{j;  a_j = \gamma, j < i \}.
\end{equation}

We define the  left ideal $\mathbf{Irr}$  of $\textbf{U}(gl(m_1|m_2+n))$ as the {\it {left ideal}} generated by the set of
irregular expressions.
\end{definition}

\begin{remark}

The action
of any element of $\mathbf{Irr}$ on the subalgebra $\mathbb C[M_{n,d}] \subset \mathbb{C}[M_{m_1|m_2+n,d}]$ - via the representation $\varrho$ -
is identically zero.
\end{remark}

\medskip

We have the  immersions
$$
\textbf{U}(gl(n)) \equiv \textbf{U}(gl(0|n)) \hookrightarrow
\textbf{U}(gl(m_1|m_2+n)),
$$
induced by the immersion
$gl(n) \hookrightarrow gl(m_1|m_2+n)$.

\medskip

\begin{lemma}

The sum ${\mathbf{U}}(gl(0|n)) + \textbf{Irr}$ is direct sum of vector spaces.
\end{lemma}

The vector space $Virt(m_1+m_2,n) = \mathbf{U}(gl(0|n)) \oplus \textbf{Irr}$
is a subalgebra of ${\mathbf{U}}(gl(m_1|m_2+n)),$ which we called the {\textit{virtual subalgebra}} \cite{Brini3-BR}.

\begin{lemma}

$\textbf{Irr}$ is a two sided ideal of $Virt(m_1+m_2,n).$

\end{lemma}

\subsection{The virtual algebra $Virt(m_1+m_2,n)$ and the virtual
presentations of elements in $\mathbf{U}(gl(n))$}

\begin{theorem}{\bf{(The Capelli epimorphism $\pi$)}}

\begin{enumerate}
\item
Every operator in $\varrho[ Virt(m_1+m_2,n) ]$ leaves invariant the
algebra of algebraic forms $\mathbb{C}[M_{n,d}] \subseteq
\mathbb{C}[M_{m_1|m_2+n,d}].$

\item
Let $n \leq d.$
The morphism
$$
Virt(m_1+m_2, n) \rightarrow
\mathrm{End}_\mathbb{C}[\mathbb{C}[M_{m_1|m_2+n,d}]] \rightarrow
\mathrm{End}_\mathbb{C}[\mathbb{C}[M_{n,d}]]
$$

defines a surjective morphism
$$
\pi : Virt(m_1+m_2, n) \rightarrow {\mathbf{P}}_{n,d} \cong
{\mathbf{U}}(gl(n)), \qquad  Ker(\pi) = \textbf{Irr}.
$$
\end{enumerate}

\end{theorem}

The ``devirtualization'' projection operator $\pi$ is called the {\it{Capelli epimorphism}}.

\begin{remark}
Any element in $Virt(m_1+m_2,n)$ defines an element in
$\mathbf{U}(gl(n))$, and  is called a \textit{virtual
presentation} of it.  The map $\pi$ being a surjection, any element
$\mathbf{p} \in \mathbf{U}(gl(n))$ admits several virtual
presentations.

From the point of view of identities, the idea of the {\it{method of virtual variables}}
may be informally summarized as follows: in order to prove
identities in $\mathbf{U}(gl(n))$
\begin{itemize}

\item

look for ``simple" virtual presentations of the elements involved in;

\item

prove the identity among the virtual presentations in
$Virt(m_1+m_2,n)$ {\it{modulo}} the ideal $\textbf{Irr}.$

\end{itemize}

\end{remark}

\begin{remark}\label{rappresentazione aggiunta-BR}

for every $e_{x_i, x_j} \in gl(n) \subset gl(m_1|m_2+n)$,   let $ad(e_{x_i, x_j})$ denote its adjoint action
on  $Virt(m_1+m_2,n)$; the ideal $\textbf{Irr}$ is $ad(e_{x_i, x_j})-$invariant. Then
\begin{equation}
\pi \left( ad(e_{x_i, x_j})( \mathbf{m} ) \right) =  ad(e_{x_i, x_j}) \left( \pi ( \mathbf{m} ) \right) ,
\qquad  \mathbf{m} \in Virt(m_1+m_2,n).
\end{equation}
\end{remark}

\begin{definition}{\bf{(Balanced monomials)}}\label{balanced monomials-BR}

In the enveloping algebra ${\mathbf{U}}(gl(m_1|m_2+n)),$ consider an
element of the form:
\begin{equation}
e_{x_{i_1},\gamma_{p_1}} \cdots e_{x_{i_n},\gamma_{p_n}} \cdot
e_{\gamma_{p_1},x_{j_1}} \cdots e_{\gamma_{p_n},x_{j_n}}, \qquad
\gamma_{p_k} \in A_0 \cup A_1,
\end{equation}
$$(x_{i_1}, \ldots, x_{i_n}, x_{j_1}, \ldots, x_{j_n} \in L, \
i.e.
proper \ symbols)
$$
that is an element that {\it creates} some virtual symbols
$\gamma_{p_1}, \ldots, \gamma_{p_n}$ (with  prescribed
multiplicities) {\it times} an element that {\it annihilates} the
{\it same} virtual symbols (with the {\it same} prescribed
multiplicities).

We call such a monomial a \textit{balanced monomial}.

\end{definition}

\begin{proposition}\footnote{This result is the (superalgebraic) formalization of the argument developed by Capelli in
\cite{Cap4-BR}, CAPITOLO I, §X.Metodo delle variabili ausiliarie, page $55$ ff.}
Every balanced monomial belongs to $Virt(m_1+m_2,n)$. Hence
$$
\mathbf{p} = \pi \left[ e_{x_{i_1},\gamma_{p_1}} \cdots
e_{x_{i_n},\gamma_{p_n}} \cdot e_{\gamma_{p_1},x_{j_1}} \cdots
e_{\gamma_{p_n},x_{j_n}} \right] \in \textbf{U}(gl(n)),
$$
Furthermore, $\mathbf{p} \in \textbf{U}(gl(n))$ is independent of the choice
of the virtual symbols $\gamma_{p_h}.$
\end{proposition}

We will say that the balanced monomial
$$e_{x_{i_1},\gamma_{p_1}} \cdots e_{x_{i_n},\gamma_{p_n}} \cdot
e_{\gamma_{p_1},x_{j_1}} \cdots e_{\gamma_{p_n},x_{j_n}}$$
is a {\it{monomial virtual presentation}} of the element
$\mathbf{p} \in \textbf{U}(gl(n)).$

\medskip

The following result lies deeper and  is a major tool in the proof of identities involving
monomial virtual presentation of elements of $\textbf{U}(gl(n)).$ Since the adjoint representation acts by superderivation,
it may be regarded as a version of the {\it{Laplace expansion}} for the images of balanced monomials.

\begin{proposition}{\bf{(Monomial virtual presentation and adjoint actions)}}\label{Monomial virtual presentation and adjoint actions-BR}

In $\textbf{U}(gl(n)),$ the element
$$
\pi \left[ e_{x_{i_1},\gamma_{p_1}} \cdots e_{x_{i_n},\gamma_{p_n}} \cdot
e_{\gamma_{p_1},x_{j_1}} \cdots e_{\gamma_{p_n},x_{j_n}} \right]
$$
equals
$$
\pi \left[ ad(e_{x_{i_1},\gamma_{p_1}}) \cdots ad(e_{x_{i_n},\gamma_{p_n}}) \left( e_{\gamma_{p_1},x_{j_1}} \cdots
e_{\gamma_{p_n},x_{j_n}} \right) \right].
$$

\end{proposition}

\begin{example}

Let $\alpha \in A_1$. The element
$$\mathbf{p} =  \pi \left[ e_{x_{3},\alpha} e_{x_{2},\alpha} e_{x_{1},\alpha} \cdot
e_{\alpha,x_{1}} e_{\alpha,x_{2}} e_{\alpha,x_{3}} \right] =
$$
$$
= \pi \left[ ad(e_{{x_3},\alpha}) ad(e_{{x_2},\alpha}) ad(e_{{x_1},\alpha}) \left( e_{\alpha,x_1} e_{\alpha,x_2} e_{\alpha,x_3} \right) \right]
$$

equals the {\it{column permanent}}\footnote{The symbol $\textbf{cper}$
denotes the column permanent of a matrix $A = [a_{ij}]$ with noncommutative entries:
$\textbf{cper} (A) = \sum_{\sigma} \ a_{\sigma(1), 1}a_{\sigma(2), 2} \cdots a_{\sigma(n), n}.$}

$$
\textbf{cper}
\left(
 \begin{array}{ccc}
 e_{{x_1},{x_1}} - 2 & e_{{x_1},{x_2}} &  e_{{x_1},{x_3}} \\
 e_{{x_2},{x_1}} & e_{{x_2},{x_2}} - 1 &  e_{{x_2},{x_3}}\\
 e_{{x_3},{x_1}} & e_{{x_3},{x_2}} &  e_{{x_3},{x_3}}\\
 \end{array}
 \right) \in {\mathbf{U}}(gl(3)).
$$
\end{example}

\section{Capelli bitableaux and the Koszul map}
Let $S$ and $T$ be two Young tableaux of same shape $\lambda = (\lambda_1 \geq \lambda_2 \geq \cdots \geq \lambda_p)$ on the same alphabet $L = \{ x_1, \ldots, x_n \}$:

\begin{equation}\label{bitableaux}
S = \left(
\begin{array}{llllllllllllll}
x_{i_1}  \ldots    \ldots     \ldots    & x_{i_{\lambda_1}}     \\
x_{j_1}   \ldots  \ldots               x_{j_{\lambda_2}} \\
 \ldots  \ldots  & \\
x_{s_1} \ldots x_{s_{\lambda_p}}
\end{array}
\right), \qquad
T = \left(
\begin{array}{llllllllllllll}
x_{h_1}  \ldots    \ldots     \ldots    & x_{h_{\lambda_1}}    \\
x_{k_1}   \ldots  \ldots               x_{k_{\lambda_2}} \\
 \ldots  \ldots  & \\
x_{t_1} \ldots x_{t_{\lambda_p}}
\end{array}
\right)
\end{equation}
and let  $\gamma_1, \ldots, \gamma_p$ be $p$ different virtual symbols of the
{\it{same parity}}.

In ${\mathbf{U}}(gl(m_1|m_2+n)),$ we consider
the monomials:
$$
e_{S, \gamma} =
\begin{array}{llllllllllllll}
e_{x_{i_1}, \gamma_1}  \cdots    \cdots     \cdots    & e_{x_{i_{\lambda_1}}, \gamma_1} \cdot    \\
e_{x_{j_1}, \gamma_2}   \cdots  \cdots               e_{x_{j_{\lambda_2}}, \gamma_2} \cdot  \\
 \cdots  \cdots  & \cdot \\
e_{x_{s_1}, \gamma_p} \cdots e_{x_{s_{\lambda_p}}, \gamma_p}
\end{array}
\
$$
and
$$
e_{\gamma, T} =
\begin{array}{llllllllllllll}
e_{\gamma_1,x_{h_1}}  \cdots    \cdots     \cdots    & e_{\gamma_1, x_{h_{\lambda_1}}} \cdot    \\
e_{\gamma_2, x_{k_1}}   \cdots  \cdots               e_{\gamma_2, x_{k_{\lambda_2}}} \cdot  \\
 \cdots  \cdots  & \cdot \\
e_{\gamma_p, x_{t_1}} \cdots e_{\gamma_p, x_{t_{\lambda_p}}}
\end{array}.
$$
and we define the balanced monomial:
\begin{equation}
e_{S, \gamma,  T} = e_{S,\gamma} \cdot e_{\gamma, T} \in
Virt(m_1+m_2,n) \subseteq {\mathbf{U}}(gl(m_1|m_2+n)).
\end{equation}

\begin{definition}{\bf{(Determinantal and permanental Capelli
bitableaux)}}

If $\gamma_1, \ldots, \gamma_p \in A_0$, the element
$$
[S|T] = \pi \left( e_{S, \gamma,  T} \right) \in {\mathbf{U}}(gl(n))
$$
is called a {\it{determinantal Capelli bitableau}}.

If $\gamma_1, \ldots, \gamma_p \in A_1$, the element
$$
[S|T]^{*} = \pi \left( e_{S, \gamma,  T} \right) \in
{\mathbf{U}}(gl(n))
$$
is called a {\it{permanental Capelli bitableau}}.

The balanced monomials $e_{S, \gamma,  T}, \ \gamma \in A_0$ ($\gamma \in A_1$) are  {\it{monomial virtual presentations}}
of the determinantal Capelli bitableau $[S|T]$ ( permanental Capelli bitableau $[S|T]^{*} )$.

\end{definition}

By referring to their monomial virtual presentations, we see that determinantal (permanental) Capelli bitableaux are
{\it{skew-symmetric}}  ({\it{symmetric}}) with respect to permutations of elements in the same row, both in the tableaux $S$ and $T.$

\begin{example}\label{Capelli determinants-BR}

Let $\alpha \in A_0$. Then
$$
 [ x_{i_k} \cdots x_{i_2}  x_{i_1} | x_{i_1} x_{i_2} \cdots x_{i_k} ] =
 \pi \left( e_{x_{i_k},\alpha} \cdots e_{x_{i_2},\alpha} e_{x_{i_1},\alpha} \cdot e_{\alpha, x_{i_1}} e_{\alpha, x_{i_2}}
 \cdots e_{\alpha, x_{i_k}} \right) =
 $$
 $$
= \textbf{cdet}\left(
 \begin{array}{cccc}
 e_{x_{i_1},x_{i_1}}+(k-1) & e_{x_{i_1},x_{i_2}} & \ldots  & e_{x_{i_1},x_{i_k}} \\
 e_{x_{i_2},x_{i_1}} & e_{x_{i_2},x_{i_2}}+(k-2) & \ldots  & e_{x_{i_2},x_{i_k}}\\
 \vdots  &    \vdots                            & \vdots &  \\
e_{x_{i_k},x_{i_1}} & e_{x_{i_k},x_{i_2}} & \ldots & e_{x_{i_k},x_{i_k}}\\
 \end{array}
 \right) \in {\mathbf{U}}(gl(n)).
$$
\end{example}

\begin{example}
Let $\alpha, \beta \in A_0$. Then
$$
\left[
\begin{array}{c|c}
x_1 & x_2 \\
x_2 & x_1
\end{array}
\right] = \pi \left( e_{x_1,\alpha}e_{x_2,\beta} \cdot e_{\alpha,x_2}e_{\beta,x_1} \right) =
$$
$$
= - e_{x_1,x_2}e_{x_2,x_1} + e_{x_1,x_1} \in {\mathbf{U}}(gl(2)).
$$
\end{example}

\subsection{The Koszul map and the bitableaux isomorphism}

Let $\mathbb{C}[ (x_i|x_j) ]_{i, j = 1, \ldots, n}$ denote the polynomial $\mathbb{C}-$algebra in the
variables $(x_i|x_j),$  $i, j = 1, \ldots, n$;  this algebra is isomorphic to  $Sym[gl(n)],$
via the map
$$
(x_i|x_j) \mapsto e_{x_{i},x_{j}}, \quad i, j = 1, \ldots, n.
$$

For every $x_h, x_k$, $h, k = 1, 2, \ldots, n$ let
$$
\rho_{x_h, x_k} : \mathbb{C}[ (x_i|x_j) ]_{i, j = 1, \ldots,n} \longrightarrow \mathbb{C}[ (x_i|x_j) ]_{i, j = 1, \ldots,n}
$$
be the linear map such that
$$
\rho_{x_h, x_k}(\mathbf{M}) = D_{x_h, x_k}(\mathbf{M}) + (x_h|x_k) \cdot \mathbf{M},
$$
for every $\mathbf{M} \in \mathbb{C}[ (x_i|x_j) ]_{i, j = 1, \ldots,n}$
( here, the symbol $D_{x_h, x_k}$ denotes the polarization  operator
$$
D_{x_h, x_k} = \sum_{j=1}^n \ (x_h|x_j) \ \frac {\partial} {\partial
(x_k|x_j)}
$$
on the algebra $\mathbb{C}[ (x_i|x_j) ]_{i, j = 1, \ldots,n}.$ )

\medskip

The map
$
e_{x_h, x_k} \rightarrow \rho_{x_h, x_k}
$
defines a Lie algebra homomorphism
$$
T: gl(n) \longrightarrow End_{\mathbb{C}}(\mathbb{C}[ (x_i|x_j) ]).
$$
By the universal property of the enveloping algebra $\mathbf{U}(gl(n)),$ the map $T$ uniquely extends to
an homomorphism of associative algebras
$$
\tau : \mathbf{U}(gl(n)) \longrightarrow \left(End_{\mathbb{C}}(\mathbb{C}[ (x_i|x_j) ]), \circ \right)
$$
such that $\tau( e_{x_h, x_k} ) = \rho_{x_h, x_k}$ and $\tau(\mathbf{1}) = Id.$

Let now  $1$ denote the unit of $\mathbb{C}[ (x_i|x_j) ]_{i, j = 1, \ldots,n}$ and let
$$
\varepsilon_1 : End_{\mathbb{C}}(\mathbb{C}[ (x_i|x_j) ]) \longrightarrow \mathbb{C}[ (x_i|x_j) ]
$$
be $\mathbb{C}-$linear map {\it{evaluation at 1}}, that is $\varepsilon_1(\rho) = \rho(1),$ for every
$\rho \in End_{\mathbb{C}}(\mathbb{C}[ (x_i|x_j) ]).$

The composite map
$$
\mathfrak{K} = \varepsilon_1 \circ \tau : \mathbf{U}(gl(n)) \longrightarrow \mathbb{C}[ (x_i|x_j) ] \cong Sym[gl(n)]
$$
is the {\it{Koszul operator}} \cite{Koszul-BR}.

We recall that, given a pair $(S, T)$ of the same shape, its {\it{biditerminant}} $(S|T)$ is an element
of $\mathbb{C}[ (x_i|x_j) ].$\footnote{In this note we assume the ``signed" definition of   \cite{rota-BR}, which
differs from those of \cite{DEP-BR}, \cite{DKR-BR}, \cite{drs-BR} just for a sign.
By referring to the notation of eq. (\ref{bitableaux}),
$$
(S|T) = \theta \cdot det[(x_{i_{\lambda_1-u+1}}|x_{h_{u'}})]_{u,u'=1, \ldots, \lambda_1} \cdots
det[(x_{s_{\lambda_p-v+1}}|x_{t_{v'}})]_{v,v'=1, \ldots, \lambda_p} \in \mathbb{C}[ (x_i|x_j) ],
$$
where
$\theta = (-1)^{\lambda_1(\lambda_2 + \cdots + \lambda_p)+ \cdots + \lambda_{p-1}\lambda_p}.$
For a simple presentation of the bideterminant and of its supersymmetric analogues in terms of virtual variables,
we refer the reader to \cite{Brini1-BR}.}

\begin{theorem} (The standard basis theorem, see e.g. \cite{DEP-BR}, \cite{DKR-BR}, \cite{drs-BR})

The set of bideterminants
$$
\{ (S|T); \ S, T semistandard \footnote{A Young tableau is said to be semistandard if its rows are left to right increasing sequences and its
columns are top to bottom nondecreasing sequences.} \}
$$
is a  linear basis of $\mathbb{C}[ (x_i|x_j) ].$
\end{theorem}

We now just state a series of  results, which come
as an application of the virtual variable method, and refer  the
reader to \cite{Brini3-BR} and \cite{Brini4-BR} for further details.

\begin{theorem} (\cite{Brini3-BR}, \cite{Brini4-BR})

The map\footnote{This map  is called the ``bitableuax correspondence" in \cite{Brini4-BR}.}
$$
\mathfrak{T}: (S|T) \mapsto [S|T]
$$
defines a linear invertible operator
$$
 \mathbb{C}[ (x_i|x_j) ]_{i, j = 1, \ldots,n} \rightarrow {\mathbf{U}}(gl(n)),
$$
whose  inverse is the the Koszul operator $\mathfrak{K}.$

Therefore, the Koszul operator $\mathfrak{K}$ is invertible.

%che risulta essere l'inverso dell'operatore di Koszul  (J.P.  Koszul, C.R. Acad. Sci. Paris (1981) )
\end{theorem}

\begin{corollary}(Koszul, \cite{Koszul-BR})

$$
\mathfrak{K} \left( [x_n \ldots x_2 x_1 | x_1 x_2 \ldots  x_n] \right) = (x_n \ldots x_2 x_1 | x_1 x_2 \ldots  x_n),
$$
where
$$
(x_n \ldots x_2 x_1 | x_1 x_2 \ldots  x_n) = det \left[ (x_i|x_j) \right]_{i, j = 1, 2, \ldots, n}
$$
\end{corollary}

Besides the notion of the bideterminant $(S|T)$ of the pair $(S, T)$, one has its natural {\it{symmetric}} counterpart, namely,
its {\it{bipermanent}}  $(S|T)^{*} \in \mathbb{C}[ (x_i|x_j) ]$ (see, e.g \cite{Brini1-BR}, \cite{Bri-BR}).

We recall that the set of co-semistandard\footnote{A Young tableau is said to be co-semistandard if its rows are left
to right nondecreasing sequences and its columns are top to bottom increasing sequences.} bipermanents
is a  linear basis of $\mathbb{C}[ (x_i|x_j) ]$ (see, e.g. \cite{rota-BR}, \cite{Bri-BR}).

\begin{theorem} (\cite{Brini3-BR}, \cite{Brini4-BR})

We have:
$$
\mathfrak{T}: (S|T)^{*} \mapsto [S|T]^{*},
$$
for every $(S|T)^{*} \in \mathbb{C}[ (x_i|x_j) ].$
\end{theorem}

\begin{corollary} (\cite{Brini3-BR}, \cite{Brini4-BR})
\begin{enumerate}
\item
The set of determinantal Capelli bitableaux:
\begin{equation}\label{eq basis-BR}
\{ [S|T]; \ S, T \ semistandard \}
\end{equation}
is a linear basis of ${\mathbf{U}}(gl(n)).$

\item
The set of permanental Capelli bitableaux:
\begin{equation}\label{eq basis-BR}
\{ [S|T]^{*}; \ S, T \ co-semistandard  \}
\end{equation}
is a linear basis of ${\mathbf{U}}(gl(n)).$

\end{enumerate}
\end{corollary}

\section{The center $\mathfrak{Z}(gl(n))$ of $\mathbf{U}(gl(n))$}

In order to make  the notation lighter, in this section we simply write $1, 2, \ldots, n$ in place of
$x_1, x_2, \ldots, x_n.$

\begin{remark}\label{centrality-BR}
Throughout this section  the role of Remark \ref{rappresentazione aggiunta-BR}  is ubiquitous: in order to prove that an element is central
in $\mathbf{U}(gl(n))$, we simply claim that its virtual presentation in $Virt(m_1+m_2,n)$ is annihilated by the adjoint actions
$ad(e_{i, j}), \  e_{i, j} \in gl(n).$
\end{remark}

\begin{example}
The prototypical example is that of the Capelli element
$$
[n \ldots 2 1|1 2 \ldots n] = \textbf{cdet} \left[  e_{i, j} + \delta_{ij}( n - i)     \right]_{i, j =1, \ldots, n} \in \mathbf{U}(gl(n)),
$$
and we discuss it in detail.

First, recall that
$ad(e_{i, j})(e_{h, \alpha}) = \delta_{j h}e_{i, \alpha}, \ ad(e_{i, j})(e_{\alpha, k}) = - \delta_{k i}e_{\alpha, j},$
for every virtual symbol $\alpha$, and that $ad(e_{i, j})$ acts as a derivation.

The monomial
$$
M = e_{n, \alpha} \cdots e_{2, \alpha} e_{1, \alpha}
e_{\alpha, 1}e_{\alpha, 2} \cdots e_{\alpha, n} \in Virt(m_1+m_2,n), \quad \alpha \in A_0,
$$
is annihilated by
$ad(e_{i, j}), \ i \neq j,$ by skew-symmetry. Furthermore, $ad(e_{i, i})(M) = M - M = 0, \ i = 1, 2, \ldots, n.$

Since
$[n \ldots 2 1|1 2 \ldots n] = \pi \left( e_{n, \alpha} \cdots e_{2, \alpha} e_{1, \alpha}
e_{\alpha, 1}e_{\alpha, 2} \cdots e_{\alpha, n} \right),
\ \ \alpha \in A_0,$ the element $[n \ldots 2 1|1 2 \ldots n]$ is central in $\mathbf{U}(gl(n)),$
by Remark \ref{rappresentazione aggiunta-BR}.
\end{example}

\subsection{The classical Capelli generators of $1893$}

In the enveloping algebra $\mathbf{U}(gl(n))$, given any integer $k = 1, 2, \ldots, n,$ consider the element
(compare with {\bf{Example}} \ref{Capelli determinants-BR})
\begin{equation}\label{The classical Capelli generators of $1893$-BR}
%diag(1-(n-i))+(D_{x_i,x_j})
\mathbf{H}_n^{(k)} =
 \sum_{1 \leq i_1 < \cdots < i_k \leq n} \ [ i_k \cdots i_2 i_1 | i_1 i_2 \cdots i_k ] =
\end{equation}
$$
 = \sum_{1 \leq i_1 < \cdots < i_k \leq n} \ \textbf{cdet}\left(
 \begin{array}{cccc}
 e_{{i_1},{i_1}}+(k-1) & e_{{i_1},{i_2}} & \ldots  & e_{{i_1},{i_k}} \\
 e_{{i_2},{i_1}} & e_{{i_2},{i_2}}+(k-2) & \ldots  & e_{{i_2},{i_k}}\\
 \vdots  &    \vdots                            & \vdots &  \\
e_{{i_k},{i_1}} & e_{{i_k},{i_2}} & \ldots & e_{{i_k},{i_k}}\\
 \end{array}
 \right).
 $$

Notice that the elements $\mathbf{H}_n^{(k)}$ are easily seen to be central, by Remark \ref{centrality-BR}.

We recall the following fundamental result, proved by  Capelli in  two  papers (\cite{Cap2-BR}, \cite{Cap3-BR}) with deceiving titles
(for a faithful description of
Capelli's original proof, quite simplified by means
of the superalgebraic method of virtual variables, see \cite{Brini4-BR}).

\begin{theorem}

The set
$$
\mathbf{H}_n^{(1)}, \mathbf{H}_n^{(2)}, \ldots, \mathbf{H}_n^{(n)} = \mathbf{H}_n
$$
is a set of algebraically independent generators of the center $\mathfrak{Z}(gl(n))$ of
$\mathbf{U}(gl(n)).$
\end{theorem}

\subsection{Rectangular Capelli/Deruyts bitableaux}
Given any positive integer $p$, we define the {\it{rectangular Capelli/Deruyts bitableau}}, with $p$ rows:
$$
\mathbf{K_n^p} =
\left[
\begin{array}{l}
n \ n-1 \ \ldots \ 3 \ 2 \ 1\\
n \ n-1 \ \ldots \ 3 \ 2 \ 1\\
\cdots\\
\\
\cdots\\
n \ n-1 \ \ldots \ 3 \ 2 \ 1\\
\end{array}\\
\right| \left.
\begin{array}{l}
1 \ 2 \ 3 \ \dots \ n-1 \ n \\
1 \ 2 \ 3 \ \dots \ n-1 \ n  \\
\cdots\\
\\
\cdots\\
1 \ 2 \ 3 \ \dots \ n-1 \ n  \\
\end{array}
\right] \in \mathbf{U}(gl(n)).
$$

 From Remark \ref{centrality-BR}, we  infer:

\begin{proposition}

The elements $\mathbf{K_n^p}$ are central in $\mathbf{U}(gl(n))$.
\end{proposition}

Set, by definition, $\mathbf{K_n^0} = \mathbf{1}.$

\medskip

Any rectangular Capelli/Deruyts bitableau $\mathbf{K_n^p}$ well behaves on highest weight vectors and therefore, being central,
on irreducible representations. The following result directly follows by iterating the first assertion of Proposition $5$ of \cite{Regonati-BR}.

\begin{proposition}\label{The hook coefficient lemma-BR}({\bf{The hook coefficient lemma}})

Let $v_{\mu}$ a highest weight vector of weight $\mu = (\mu_1 \geq \mu_2 \geq \ldots \geq \mu_n),$ with $\mu_i \in \mathbb{N}$
for every $i = 1, 2, \ldots, n.$
Then
$$
\mathbf{K_n^p}(v_{\mu}) = (-1)^{{p \choose 2}n} \ \left( \prod_{i = 0}^{p -1} \ (\mu_1 - i + n -1)(\mu_2 - i + n -2) \cdots
(\mu_n - i) \right)  \cdot v_{\mu}.
$$

\end{proposition}

The crucial result in this section is that rectangular Capelli/Deruyts bitableaux $\mathbf{K_n^p}$  expand - in
a beautiful way - in terms of  the classical Capelli generators     \label{The classical Capelli generators of $1893$-BR}.

\begin{theorem}({\bf{Expansion Theorem}}) \label{Teorema di sfilamento-BR}\footnote{This is a new result. It is a special case of a
 more general result - joint work with F. Regonati and
A. Teolis - that will appear  in a forthcoming publication.}

Let $p \in \mathbb{N}$ and set $\mathbf{H}_n^{(0)} = \mathbf{1},$ by definition. The following identity in $\mathfrak{Z}(gl(n))$
holds:
\begin{equation} \label{equazione di sfilamento-BR}
 \mathbf{K_n^{p + 1}} =  (-1)^{np} \ \mathbf{K_n^p} \ \mathbf{C}_n(p),
\end{equation}
where
\begin{equation} \label{relazioni lineari-BR}
\mathbf{C}_n(p) =  \sum_{j = 0}^n \ (-1)^{n - j} (p)_{n - j} \ \mathbf{H}_n^{(j)},
\end{equation}
and
$$
(p)_k = p(p-1) \cdots (p-k+1), \ p, k \in \mathbb{N}
$$
denotes the falling factorial coefficient.
\end{theorem}

For $p = 0$, the preceding  identity consistently collapses to
$$
\mathbf{K_n^{1}} = \mathbf{H}_n^{(n)} = \mathbf{H}_n = \mathbf{C}_n(0).
$$

Furthermore, notice that the linear relations (\ref{relazioni lineari-BR}), for $p = 0, \ldots, n-1$, yield a nonsingular triangular
coefficients matrix.

\begin{corollary}

The set
$
\mathbf{C}_n(0), \mathbf{C}_n(1), \ldots, \mathbf{C}_n(n - 1)
$
is a set of algebraically independent generators of the center $\mathfrak{Z}(gl(n))$ of
$\mathbf{U}(gl(n)).$
\end{corollary}

\subsection{The Harish-Chandra isomorphism and the algebra $\Lambda^*(n)$ of shifted symmetric polynomials}

In this subsection we follow  A. Okounkov and G. Olshanski
\cite{OkOlsh-BR}.

As in the classical context of the algebra $\Lambda(n)$ of symmetric
polynomials in $n$ variables $x_1, x_2, \ldots, x_n$, the algebra
$\Lambda^*(n)$ of {\it{shifted symmetric polynomials}} is an algebra
of polynomials $p(x_1, x_2, \ldots, x_n)$  but the ordinary symmetry
is replaced by the {\it{shifted symmetry}}:
$$
 f(x_1, \ldots , x_i, x_{i+1}, \ldots, x_n) = f(x_1, \ldots , x_{i+1} - 1, x_i + 1,
 \ldots, x_n),
$$
for $i = 1, 2, \ldots, n - 1.$

\vskip 0.3cm

\begin{examples}\label{Ex shifted-BR}

Two basic classes of shifted symmetric polynomials are provided by the sequences of {\it{shifted elementary symmetric polynomials}} and
{\it{shifted complete symmetric polynomials}}.

\begin{itemize}

\item {\bf{Shifted elementary symmetric polynomials}}

For every $r \in \mathbb{N}$ let
\begin{equation}\label{shifted elementary-BR}
\mathbf{e}_r^{*}(x_1, x_2, \ldots, x_n)
= \sum_{1 \leq i_1 < i_2 < \cdots < i_r \leq n} \ (x_{i_1}  + r  - 1)
(x_{i_2}  + r - 2) \cdots (x_{i_r}),
\end{equation}
and $\mathbf{e}_0^{*}(x_1, x_2, \ldots, x_n) = \mathbf{1}.$

\item {\bf{Shifted complete  symmetric polynomials}}

For every $r \in \mathbb{N}$ let
\begin{equation}\label{shifted complete-BR}
\mathbf{h}_r^{*}(x_1, x_2, \ldots, x_n)
= \sum_{1 \leq i_1 \leq i_2 < \cdots \leq i_r \leq n} \ (x_{i_1}  - r  + 1)
(x_{i_2}  - r + 2) \cdots (x_{i_r}),
\end{equation}
and $\mathbf{h}_0^{*}(x_1, x_2, \ldots, x_n) = \mathbf{1}.$

\end{itemize}

\end{examples}

\vskip 0.3cm

The {\it{Harish-Chandra isomorphism}} is the algebra isomorphism
$$
\chi : \mathfrak{Z}(gl(n)) \longrightarrow \Lambda^*(n), \qquad  \ A \mapsto \chi(A),
$$
$\chi(A)$ being the shifted symmetric polynomial such that, for every highest weight module $V_{\mu}$,
the evaluation $\chi(A)(\mu_1, \mu_2, \ldots , \mu_n)$ equals the eigenvalue of
$A \in \mathfrak{Z}(gl(n))$ in  $V_{\mu}$ (\cite{OkOlsh-BR}, Proposition $\mathbf{2.1}$).

\subsection{The Harish-Chandra isomorphism interpretation of Proposition \ref{The hook coefficient lemma-BR} and
Theorem \ref{Teorema di sfilamento-BR}}

Notice that
$$
\chi(\mathbf{H}_n^{(r)}) =  \mathbf{e}_r^{*}(x_1, x_2, \ldots, x_n) \in \Lambda^*(n),
$$
for every $r =  1, 2, \ldots, n.$

Furthermore, from Proposition \ref{The hook coefficient lemma-BR}  it  follows

\begin{corollary}

We have:
\begin{equation}\label{eq quozienti-BR}
\chi(\mathbf{K_n^p}) = (-1)^{{p \choose 2}n} \left(  \prod_{i = 0}^{p -1} \
(x_1 - i + n -1)(x_2 - i + n -2) \cdots (x_n - i) \right).
\end{equation}

\end{corollary}

By combining  equations (\ref{equazione di sfilamento-BR}) and (\ref{eq quozienti-BR}),
we infer:

\begin{proposition}
We have:

\begin{itemize}
\item

For every $p \in \mathbb{N}$,
\begin{equation}\label{eq fattorizzazione-BR}
\chi(\mathbf{C}_n(p)) = (x_1 - p + n - 1)(x_2 - p + n - 2) \cdots (x_n - p).
\end{equation}

\item

The set
$$
\chi(\mathbf{C}_n(0)), \ \chi(\mathbf{C}_n(1)), \ \ldots \ , \ \chi(\mathbf{C}_n(n - 1))
$$
is a system of algebraically independent generators of the
ring $\Lambda^*(n)$ of shifted symmetric polynomials in the variables $x_1, x_2, \ldots, x_n.$

\end{itemize}
\end{proposition}

We recall a standard result (for an elementary proof see e.g.  \cite{Umeda-BR}):

\begin{proposition}

For every $p \in \mathbb{N}$, the element
\begin{equation}\label{centrality Umeda-BR}
\mathbf{H}_n(p) = \textbf{cdet}\left(
 \begin{array}{cccc}
 e_{1,1} - p + (n-1) & e_{1,2} & \ldots  & e_{1,n} \\
 e_{2,1} & e_{2,2} - p + (n-2) & \ldots  & e_{2,n}\\
 \vdots  &    \vdots                            & \vdots &  \\
e_{n,1} & e_{n,2} & \ldots & e_{n,n} - p\\
 \end{array}
 \right)  =
\end{equation}
$$
= \textbf{cdet} \left[  e_{i, j} + \delta_{ij}(- p  + n - i)     \right]_{i, j =1, \ldots, n} \in \mathbf{U}(gl(n)).
$$
 is  central.
In symbols, $\mathbf{H}_n(p)  \in \mathfrak{Z}(gl(n)).$

\end{proposition}

\medskip

Equation (\ref{eq fattorizzazione-BR}) implies
$$
\chi(\mathbf{H}_n(p)) = (x_1 - p + n - 1)(x_2 - p + n - 2) \cdots (x_n - p) = \chi(\mathbf{C}_n(p)),
$$
and, therefore, the following

\begin{corollary}

For every $p \in \mathbb{N}$, we have
$$
\mathbf{H}_n(p) =   \mathbf{C}_n(p) = \sum_{j = 0}^n \ (-1)^{n - j} (p)_{n - j} \ \mathbf{H}_n^{(j)}.
$$
\end{corollary}

\medskip

Let $t$ be a variable and consider the polynomial
$$
\mathbf{H}_n(t) = \textbf{cdet}\left(
 \begin{array}{cccc}
 e_{1,1} - t + (n-1) & e_{1,2} & \ldots  & e_{1,n} \\
 e_{2,1} & e_{2,2} - t + (n-2) & \ldots  & e_{2,n}\\
 \vdots  &    \vdots                            & \vdots &  \\
e_{n,1} & e_{n,2} & \ldots & e_{n,n} - t\\
 \end{array}
 \right) =
$$
$$
= \textbf{cdet} \left[  e_{i, j} + \delta_{ij}(- t  + n - i)     \right]_{i, j =1, \ldots, n}
$$
with coefficients in $\mathbf{U}(gl(n)).$

\medskip

\begin{corollary}(see, e.g. \cite{Umeda-BR})\label{eqUMEDA-BR}

We have the generating function identity:
$$
 \mathbf{H}_n(t)  = \sum_{j = 0}^n \ (-1)^{n - j} (t)_{n - j} \ \mathbf{H}_n^{(j)},
$$
where, for every $k \in \mathbb{N}$, $(t)_k = t(t-1) \cdots (t-k+1)$ denotes the $k-$th falling factorial polynomial.
\end{corollary}

\begin{corollary}
We have the generating function identity:

$$
\sum_{j = 0}^n \ (-1)^{n - j} (t)_{n - j} \ \mathbf{e}_j^{*}(x_1, x_2, \ldots, x_n) =
(x_1 - t + n - 1)(x_2 - t + n - 2) \cdots (x_n - t).
$$
\end{corollary}

\medskip

Following Molev \cite{Molev-BR} Chapt. {\bf 7} (see also   Brundan \cite{Brundan-BR}, Howe and Umeda \cite{HU-BR}), consider the ``Capelli determinant''
$$
{\mathcal C}_n(s) = \textbf{cdet}\left(
 \begin{array}{cccc}
 e_{1,1} + s  & e_{1,2} & \ldots  & e_{1,n} \\
 e_{2,1} & e_{2,2} + s - 1 & \ldots  & e_{2,n}\\
 \vdots  &    \vdots                            & \vdots &  \\
e_{n,1} & e_{n,2} & \ldots & e_{n,n} + s - (n-1)\\
 \end{array}
 \right) =
$$
$$
= \textbf{cdet} \left[  e_{i, j} + \delta_{ij}(s  - i + 1)     \right]_{i, j =1, \ldots, n},
$$
 regarded as a  polynomial in the variable $s$.

\medskip

By the formal (column) Laplace rule, the coefficients ${\mathcal C}_n^{(h)} \in \mathbf{U}(gl(n))$ in the expansion
$$
{\mathcal C}_n(s)  = s^n +  {\mathcal C}_n^{(1)} s^{n - 1} + {\mathcal C}_n^{(2)} s^{n - 2} + \ldots + {\mathcal C}_n^{(n)},
$$
are the sums of  the minors:
$$
{\mathcal C}_n^{(h)} = \sum_{1 \leq i_1 < i_2 < \ldots < i_h \leq n} \ {\mathcal M}_{ i_1, i_2, \ldots, i_h},
$$
where  ${\mathcal M}_{ i_1, i_2, \ldots, i_h}$ denotes the column determinant of the submatrix of the matrix ${\mathcal C}_n(0)$
obtained by selecting the rows and the columns with indices $i_1 < i_2 < \ldots < i_h.$

\medskip

Since ${\mathcal C}_n(s) = \mathbf{H}_n(-s + (n-1))$,
from  Proposition \ref{eqUMEDA-BR}
it follows:

\begin{corollary} \label{expansion-BR}

$$
{\mathcal C}_n(s)  =   \sum_{j = 0}^n \ (-1)^{n - j} (-s + (n-1))_{n - j} \ \mathbf{H}_n^{(j)}.
$$

\end{corollary}

\begin{remark} \label{Stirling-BR}

Recall that  the  polynomial sequence of {\it{falling factorials}} $\left( (x)_n \right)_{n \in \mathbb{N}}$ is
a {\it{polynomial sequence of binomial type}} and that one has the basic linear relations
$$
(x)_n = \sum_{k = 0}^n \ {\mathrm{s}}(n,k) \ x^k, \qquad n \in \mathbb{N},
$$
where the symbol ${\mathrm{s}}(n,k)$ denotes the {\it{Stirling numbers of the first kind}} (see, e.g. \cite{Aigner-BR}).

A straightforward computation leads to the   linear relations:

$$
{\mathcal C}_n^{(n -h)} = (-1)^h \ \sum_{j = 0}^n \ c(n - h, j) \  (n - j)! \ \mathbf{H}_n^{(j)},
$$
where
\begin{equation}\label{coefficients-BR}
c(n - h, j) =
 \sum_{p = 0}^{n - j} \ {n - 1 \choose j + p - 1} \frac {{\mathrm{s}}(p,h)} {p!}.
\end{equation}

\end{remark}

\medskip

Since the coefficients in (\ref{coefficients-BR}) yield a nonsingular triangular matrix, Remark  \ref{Stirling-BR} implies:

\begin{corollary}

\begin{itemize}
We have:

\item
The  elements ${\mathcal C}_n^{(h)}, \ h = 1, 2, \ldots, n$ are central and provide a system of algebraically
independent generators of $\mathfrak{Z}(gl(n)).$

\item
$$
\chi({\mathcal C}_n^{(h)}) = \bar{\mathbf{e}}_h(x_1, x_2, \ldots, x_n) = \mathbf{e}_h(x_1 , x_2  - 1, \ldots, x_n - (n - 1)),
$$
where $\mathbf{e}_h$ denotes the $h-$th elementary symmetric polynomial.
\end{itemize}

\end{corollary}

\subsection{Permanental  generators}

We end this section by describing, in the virtual variables notation, the set of the preimages in $\mathfrak{Z}(gl(n))$
- with respect to the Harish-Chandra isomorphism - of the sequence of {\it{shifted complete  symmetric polynomials}}
$\mathbf{h}_r^{*}(x_1, x_2, \ldots, x_n).$

\begin{definition}

For every $r \in \mathbb{Z}^+$, set
\begin{equation}\label{perm gen-BR}
\mathcal{P}_n^{(r)} = \sum_{(i_1, i_2, \ldots, i_n) } \ (i_1! i_2! \cdots i_n!)^{-1} \ [n^{i_n} \cdots 2^{i_2} 1^{i_1} | 1^{i_1} 2^{i_2} \cdots n^{i_n} ]^{*},
\end{equation}
where the sum is extended to all $n-$tuples $(i_1, i_2, \ldots, i_n)$ such that $i_1 + i_2 + \cdots + i_n = r$ and any
$$
[n^{i_n} \cdots 2^{i_2} 1^{i_1} | 1^{i_1} 2^{i_2} \cdots n^{i_n} ]^{*}
$$
is a {\it{permanental}} Capelli bitableau with one row.
\end{definition}

The proofs of the following results are almost trivial, as a consequence of the definition of the elements $\mathcal{P}_n^{(r)}$
in terms of their {\it{monomial virtual presentations}} (equation (\ref{perm gen-BR})).

\begin{theorem}

We have:

\begin{itemize}

\item

For every $r \in \mathbb{Z}^+$,
the element $\mathcal{P}_n^{(r)}$ is central. In symbols, $\mathbf{P}_n^{(r)} \in \mathfrak{Z}(gl(n)).$

\item

For every $r \in \mathbb{Z}^+$,  $\chi( \mathcal{P}_n^{(r)} ) = \mathbf{h}_r^{*}(x_1, x_2, \ldots, x_n).$
\end{itemize}

\end{theorem}

Amazingly, the generators $\mathcal{P}_n^{(r)}$ coincide with the ``permanental generators'' of
$\mathfrak{Z}(gl(n))$ first discovered and studied - through a rather heavy machinery - by Umeda and  Hirai \cite{UmedaHirai-BR}
(see also Turnbull  \cite{TURN-BR}).

\end{document}